\documentclass[final]{siamltex}

\usepackage{epsfig,amsmath,amssymb}

\newtheorem{remark}[theorem]{Remark}

\def\Sp1{ \text{Sp}_1 }

\def\wg{\widehat g}
\def\wC{\widehat C}
\def\wV{\widehat V}

\def\bu{$\bullet$\ \ \ }

\title{
On the $\lambda$-equations for matching control laws
\thanks{This work was partially supported 
by grant CMC 9813182 from the National Science Foundation.}
}

\author{
David Auckly \thanks{Department of Mathematics, 
Kansas State University, Manhattan, KS 66506 ({\tt dav@math.ksu.edu})}  
\and Lev Kapitanski \thanks{Department of Mathematics, 
Kansas State University, Manhattan, KS 66506
 ({\tt levkapit@math.ksu.edu }) 
}}

\begin{document}

\maketitle

\begin{abstract} We discuss matching control laws for underactuated 
systems. We previously showed that this class of matching control 
laws is completely charactarized by  a linear system of 
first order partial differential equations for one set of variables 
($\,\lambda$)  
followed by a linear system of first order PDEs for the second set 
of variables ($\,\wg$, $\,\wV$). 
Here we derive a new first order system of partial differential 
equations that encodes all compatibility conditions for the 
$\,\lambda$-equations. We give four examples illustrating different 
features of matching control laws. The last example is a system with 
two unactuated degrees of freedom that admits only basic solutions 
to the matching equations. 
There are systems with  
many matching control laws where only basic solutions are 
potentially useful. We introduce a rank condition indicating 
when this is likely to be the case. 
\end{abstract}

\begin{keywords} 
nonlinear control, matching control laws, $\lambda$-equations, 
stabilization 
\end{keywords}

\begin{AMS}
93C10, 93D15
\end{AMS}

\pagestyle{myheadings}
\thispagestyle{plain}
\markboth{D. AUCKLY AND L. KAPITANSKI}{$\lambda$-EQUATIONS}

\section{Introduction }

Effective procedures for designing control laws are very important 
in nonlinear control theory. Explicit analytic formulae for control laws play a 
role similar 
to explicit solutions to differential equations. Such formulae exist 
in only a few special cases, but those that exist serve as 
simple models 
to develop and test more general techniques. 

In this paper we discuss a class of 
full state feedback control laws for underactuated  systems. 
In \cite{AKW} we showed that this class of matching control 
laws is completely charactarized by  a linear system of 
first order partial differential equations for one set of variables 
($\,\lambda$)  
followed by a linear system of first order PDEs for the second set 
of variables ($\,\wg$, $\,\wV$). These equations always have a 
simple family of solutions which we call basic solutions. 
The system of equations for the first set of variables 
($\,\lambda$-equations) is overdetermined. 
Here we derive a new first order system of partial differential 
equations that encodes all compatibility conditions for the 
$\,\lambda$-equations (we call these the $\,\nu$-equations). 
If only one degree of 
freedom is unactuated, the solutions to all these systems of PDEs  
can be completely analyzed. It is often possible to get explicit 
formulae for the solutions to these equations.  We also provide an 
example of a system with two unactuated degrees of freedom 
that has only basic solutions. There are systems with  
many matching control laws where only basic solutions are 
potentially useful. We write down a rank condition indicating 
when this is likely to be the case. 

During the last few years 
several researchers have investigated control laws 
in which the closed loop system 
assumes a certain structure. 
Numerous papers have been written on this subject, see [1] - [13] 
and the references therein. The control laws that 
form the subject of this present paper are described by equations 
(\ref{control}) and (\ref{match}). These equations were independently derived 
in \cite{H1} and \cite{AKW}. 
The $\lambda$-equations were first introduced
in paper \cite{AKW}. 
Even though the initial matching equations of \cite{H1} and 
\cite{AKW} form a highly nonlinear system of PDEs, 
introduction of the  $\lambda$ variables triangulates the system.
The system is triangulated in the sense that all solutions are obtained 
by first solving  first order linear equations for $\lambda$ 
and then 
solving first order linear equations for the remaining variables. 

This paper is organized as follows. Section 2  
reviews   matching control laws and the $\,\lambda$-equations, and 
introduces the $\,\nu$-equations. Section 3 specializes to systems with 
one unactuated degree of freedom. Sections 4, 5 and 6 contain examples 
illustrating three different features of matching control laws. 
The rank condition appears at the end of Section 5. 
Later we apply it in Section 6. In 
Section 7 we describe the final example of a system with two unactuated 
degrees of freedom. We show that this system has only basic matching 
control laws.

\section{ Matching equations}

We use the following notation.
\begin{description}
\item[\bu] $n$ is the number of the degrees of freedom of the mechanical 
system
\item[\bu] $x=(x^1,\dots,x^n)$ are configuration variables 
denoting the position of 
the  system, and $\dot x=(\dot x^1,\dots,\dot x^n)$ are the corresponding 
velocities
\item[\bu] $\,g_{ij}(x)\,$ is the mass matrix
\item[\bu] $\,V(x)\,$ is the potential energy
\item[\bu] $C_i(x,\dot x)\,$ are 
the dissipation terms
\item[\bu] $u_i(x, \dot x)$ are the control inputs
\end{description}
Let $m\le n$ be the number of unactuated degrees of freedom. We will 
assume 
that degrees of freedom numbered $1$ through $m$ are unactuated and use 
indices $a$, $b$, \dots to indicate unactuated degrees of freedom. 
The indices 
$i$, $j$, \dots will run from $1$ to $n$. 
We adopt the convention of summation over the repeated indices . 

Given this, the equations of motion of the system are
\begin{equation}
 g_{rj}\ddot x^j\,+\, [j\,k,\,r]\,\dot x^j\,\dot x^k\,+\,
 C_{r}\,+\,{\partial V\over\partial x^r}\,=\,u_r\,,
\quad r=1,\dots, n,
\label{eq}
\end{equation}
where $[i\,j,\,k]$ is the Christoffel symbol of the first kind, 
\begin{equation}
[i\,j,\,k]\,=\,\frac12\,\big(\, {\partial g_{jk}\over \partial x^i}\,+\,
{\partial g_{ik}\over \partial x^j}\,
-\,{\partial g_{ij}\over \partial x^k}\,\big)
\label{christ}
\end{equation}
Our assumtion that  the first $m$ degrees of freedom 
are not actuated means that
\begin{equation}
u_1=\dots = u_{{}_{m}}\,=\,0.
\label{0}
\end{equation}
We are looking for control laws $u_i$ such that the closed loop 
system can be written in the form 
$$
 \wg_{rj}\ddot x^j\,
+\, \widehat {[j\,k,\,r]}\,\dot x^j\,\dot x^k\,+\,
 \wC_{r}\,+\,{\partial  \wV\over\partial x^r}\,=\,0\,,
\quad r=1,\dots, n,
$$
where $\widehat{[i\,j,\,k]}$ is defined as in (\ref{christ}) 
with $\wg\,$ in place of $g$.
Such a control law will be given by 

\begin{eqnarray}
u_\ell\,=\,
\big(\,[j\,k,\,\ell]\,-\,
 g_{\ell i}\,\wg^{ir}\,
\widehat {[j\,k,\,r]}\,\big)\,
\dot x^j\,\dot x^k\, +\, 
\, \big(\, C_{\ell}\,-\,\,g_{\ell i}\,
\wg^{ij}\,\widehat C_{j}\,\big)\,\nonumber\\  
+\,\big(\,{\partial V\over\partial x^\ell}\,-\,g_{\ell i}\,
\wg^{ij}\,{\partial\widehat V \over\partial x^j}\,\big)\,,
\qquad\qquad  \ell =1,\dots, n,
\label{control}
\end{eqnarray}
Condition (\ref{0}) translates into 
\begin{eqnarray}
\big(\,[j\,k,\,a]\,-\,
 g_{ai}\,\wg^{ir}\,
\widehat {[j\,k,\,r]}\,\big)\,
\dot x^j\,\dot x^k\, +\, 
\, \big(\, C_{a}\,-\,g_{ai}\,
\wg^{ij}\,\widehat C_{j}\,\big)\, \nonumber \\  
+\,\big(\,{\partial V\over\partial x^a}\,-\,g_{ai}\,
\wg^{ij}\,{\partial\widehat V \over\partial x^j}\,\big)\,=\,0\,,
\qquad\qquad  a=1,\dots, m,
\label{controlmatch}
\end{eqnarray}
In order to satisfy these equations it is sufficient to have
\begin{eqnarray}
g_{ai}\,\wg^{ir}\;
\widehat {[j\,k,\,r]}\;& = & \;[j\,k,\,a]  \nonumber \\
g_{ai}\,\wg^{ir} \;
\widehat C_{r}\;&  = & \; C_{a}  \label{match}\\
g_{ai}\,\wg^{ir}\;
{\partial\widehat V \over\partial x^r}\; &  = & \;
{\partial V\over\partial x^a} \nonumber  
\end{eqnarray}
These are the {\it matching equations\/}, see \cite{AKW}, \cite{H1}.
Following \cite{AKW}, introduce  variables $\lambda_a^j$
relating the unknown mass matrix $\wg$  to the original mass matrix $g$,
\begin{equation}
\lambda_a^r\,=\,g_{ai}\;\wg^{ir}\,.
\label{lamdef}
\end{equation}
Using $\lambda_a^j$,\  the matching equations take the form
\begin{eqnarray}
\lambda_a^r\;
\widehat {[j\,k,\,r]}\;& = & \;[j\,k,\,a] \nonumber \\
\lambda_a^j \;
\widehat C_{j}\;&  = & \; C_{a} \label{lammatch}\\
\lambda_a^j\;
{\partial\widehat V \over\partial x^j}\; &  = & \;
{\partial V\over\partial x^a} \nonumber  
\end{eqnarray} 
\begin{theorem}\label{one}  The following equations are 
equivalent  to the matching equations .
\begin{description}
  \item[$\lambda$-equations:]
\begin{equation}
{\partial\hfil\over\partial x^k}\,(\,g_{ai}\;\lambda^i_b)\,
-\;[k\,a,\,i]\;\lambda^i_b\,
-\;[k\,b,\,i]\;\lambda^i_a\,=\,0\,,\;
\begin{array}{ll} \, & k=1,\dots, n \\
\, & a,b \,=\,1,\dots,m
\end{array}
\label{lameq}
\end{equation}
  \item[$\wg$-equations:]
\begin{equation}
\lambda^\ell_a\,{\partial \wg_{ij}  \over\partial x^\ell}  \,+\,
{\partial \lambda^\ell_a \over\partial x^i}\cdot\wg_{\ell j}\,+\,
{\partial \lambda^\ell_a  \over\partial x^j}\cdot\wg_{\ell i}\,=\,
{\partial g_{ij} \over\partial x^a} \,, 
  \begin{array}{ll}\quad  & a=1,\dots, m \\
 & i,j=1,\dots, n\,
  \end{array}
\label{gheq}
\end{equation}
 \item[$\wV$-equations:]
\begin{equation}
\lambda_a^j\;
{\partial\widehat V \over\partial x^j}\;= \;{\partial V\over\partial x^a} 
\label{vheq}
\end{equation}
 \item[$\wC$-equations:]
\begin{equation}
\lambda_a^j \;
\widehat C_{j}\;  =  \; C_{a}
\label{cheq}
\end{equation}
\end{description}

\end{theorem}

For the proof see \cite{AKW}, \cite{AK}.

\begin{remark}\label{remarkone} These equations always have a set of 
solutions of the form
\begin{equation}
\lambda^k_a\,=\,\varkappa\,\delta^k_a,\qquad\qquad\widehat g \,
= \,\frac1{\varkappa}\,g + g^o,
\qquad\qquad\widehat V \,= \,\frac1{\varkappa}\,V+V^o,\qquad \wC_j\,
=\,\frac1{\varkappa}\,C_j
\nonumber
\end{equation}
with $\varkappa\ne 0$ any constant,  $V^o$ arbitrary function 
of the variables $\,x^\ell $, $\,\ell=m+1,\dots, n$, and $g^o$ any 
symmetric 
matrix 
valued function 
of the variables $\,x^\ell $ such that  $\,g^o_{ia}=0$.
We will call these solutions {\em basic\/}.
\end{remark}

The  $\lambda$-equations are a system of 
$\,\frac{1}{2}\,m(m+1)\cdot n$ equations for $n\cdot m$ unknowns. 
It is not surprising that there are extra compatibility conditions. 
By viewing  system (\ref{lameq}) in the correct way, we are able to write 
down the compatibility conditions. 
Denote
\begin{equation}
\nu_{ab}\,=\,g_{ai}\;\lambda^i_b\,.
\label{nueq1}
\end{equation}
Because the matrix $g_{ij}$ is assumed to be non-degenerate, 
the matrix comprised of its $m$  first rows has rank $m$. 
This implies that $\,m^2\,$ out of $\,m\cdot n\,$ $\,\lambda$'s 
can be expressed as linear combinations of $\,\nu$'s, i.e.,
$$
\lambda^\beta_b\,=\,h^{\beta a}\,\nu_{ab}\,.
$$
Substituting this 
in the $\lambda$-equations, we obtain
\begin{equation}
\partial_k\,\nu_{ab}\,-\,[a\,k,\,\beta]\,h^{\beta d}\,\nu_{db}\,
-\,[b\,k,\,\beta]\,h^{\beta d}\,\nu_{da}\,=\,
[a\,k,\,\rho]\,\lambda^\rho_b\,+\,[b\,k,\,\rho]\,\lambda^\rho_a\,,
\label{nueq2}
\end{equation}
where index $\,\rho\,$ varies over the remaining $\,(n-m)\,$ indices.
We will view  system (\ref{nueq2}) of $\,\frac{1}{2}\,m(m+1)\cdot n\,$  
equations 
as a linear  algebraic system for the $\,m\,(n-m)\,$ variables 
$\,\lambda^\rho_a\,$,
\begin{equation}
\left[A_{(k,a,b)}\right]^c_\rho\;\lambda^\rho_c\,=\,F_{(k,a,b)}\,.
\label{nualg}
\end{equation}
We know that this system has at least one solution by 
Remark \ref{remarkone}. 
Thus, the rank of the matrix $\,A\,$ is at most $\,m\cdot n$. 
In order for system (\ref{nualg}) to have a solution, the vector 
$\,F_{(k,a,b)}\,$ 
must be perpendicular to the kernel of the transposed matrix, $\,A^\star$, 
\begin{equation}
F\,\perp\,\ker\,A^\star\,.
\label{fredh}
\end{equation}
Let the kernel of the matrix $\,A^\star \,$ be generated by the 
vectors $\,\xi_r\,$. 
The orthogonality condition (\ref{fredh}) then takes the form of 
the following 
system of linear first order partial differential equations for 
$\,\nu_{ab}$:
\begin{equation}
\xi_r^{(k,a,b)}(x)\,\left[\,
\partial_k\,\nu_{ab}\,-\,[a\,k,\,\beta]\,h^{\beta d}\,\nu_{db}\,
-\,[b\,k,\,\beta]\,h^{\beta d}\,\nu_{da}\,\right]\,=\;0\,.
\label{fredh2}
\end{equation}

\begin{theorem}\label{la0} The general solution 
to the $\lambda$-equations is given by any set of $\,\lambda^\rho_a$ solving the algebraic
system (\ref{nualg}), and $\,\lambda^\beta_b\,=\,h^{\beta a}\,\nu_{ab}\,$, 
where 
$\,\nu_{ab}\,$ is any solution to equations (\ref{fredh2}). 
\end{theorem}

In general, if $\,m>1$,  system (\ref{fredh}) may be quite complicated 
and we do not have a
satisfactory  description of its solutions.

\section{Systems with one unactuated degree of freedom}

If only one degree of freedom is unactuated, we do have a reasonable 
description of all solutions to system (\ref{fredh2}). 
Assume, for simplicity, that $\,g_{11}(x) > 0$.  
Then, after rescaling $x^1$ if necessary, we will have $\,g_{11}(x)\,=\,1$. 
More precisely, from the very beginning we could use, instead of 
$\,(x^1, x^2, \dots, x^n)$, the coordinates $\,(z^1, z^2, \dots , z^n)\,$ 
which are related to $x$ as follows:
$$
{\partial z^1\over \partial x^1}\,=\,\sqrt{g_{11}(x)},\qquad 
z^2\,=\,x^2,\quad\dots ,\;z^n\,=\,x^n
$$
In $z$ coordinates the mass matrix is 
$\,\tilde g_{ij}(z)\,=\,g_{k\ell}(x)\,{\partial x^k\over\partial z^i}\,
{\partial x^\ell\over\partial z^j}$ and, hence, $\tilde g_{11}(z)=1$. 
On the other hand, the structure of the equations of motion 
(\ref{eq}) does not change because of their tensorial form, and 
the condition $u_1=0$ remains the same, again, because 
$\,\tilde u_1\,=\,u_k\,{\partial x^k\over\partial z^1}\,
=\,u_1\,\sqrt{g_{11}(x)}$. 
Thus, we assume that the coordinates are chosen appropriately, and 
$\,g_{11}(x)=1$. 

In the case of one unactuated degree of freedom one is solving for 
$\,\lambda^i_1$. The $\,\lambda$-equation reads
\begin{equation}
{\partial\nu\over\partial x^k}\,=\,2\,[k\,1,\,i]\,\lambda^i_1\,,
\label{nu}
\end{equation}
where $\,\nu\,=\,g_{1i}\,\lambda^i_1$. 
Notice that $\,[k\,1,\,1]\,=\,0$. View the equations (\ref{nu}) 
as a system of linear algebraic equations for the variables 
$\,\lambda^\rho_1$, $\rho=2,\dots, n$. In order for this system of 
$n$ equations in $(n-1)$ unknowns to have a solution,  the vector 
$$
\,v\,=\,
\begin{pmatrix} 
\partial_1\,\nu \\ \dots \\ \partial_n\nu 
\end{pmatrix}
$$ 
must be 
perpendicular to the kernel of the matrix
\begin{equation}
A^\star\,=\,
\begin{pmatrix} 
\,[1\,1,\,2] & \dots & \,[n\,1,\,2] \\ 
  \dots & \dots & \dots \\ 
\,[1\,1,\,n] & \dots & \,[n\,1,\,n]
\end{pmatrix}\qquad .
\nonumber
\end{equation} 
Let the kernel of $\,A^\star\,$ be generated by the vectors 
$\,\xi_r\,=\,(\xi_r^1,\dots,\,\xi_r^n)$. The orthogonality 
condition for $v$ translates into the system of equations
\begin{equation}
X_r(\nu)\,\equiv\,
\xi_r^1(x)\,{\partial\nu\over\partial x^1}\,+\,\dots\,+
\,\xi_r^n(x)\,{\partial\nu\over\partial x^n}\,=\,0.
\label{nu2}
\end{equation}
The standard procedure to solve such a system of equations is to 
complete the system into an involutive system by  
adding equations 
$\,[X_r,\,X_s](\nu)\,=\,0$, $\,[[X_r,\,X_s],\,X_t](\nu)\,=\,0$, \dots,
where 
$\,[\eta^i\,\partial_i,\,\zeta^j\,\partial_j]\,
=\,\big(\eta^i\partial_i(\zeta^k)\,
-\,\zeta^i\,\partial_i(\eta^k)\big)\,\partial_k$ is the commutator of 
vector-fields. 
Recall that a system of equations 
$$
Y_1(\nu)\,=\,0,\dots,\,Y_K(\nu)\,=\,0
$$
is involutive if $\,[Y_p,\,Y_q]\,=\,f^r_{pq}(x)\,Y_r$.
 
Thus we have proved the following result.
\begin{theorem}\label{la1} With one 
unactuated degree of freedom there is a coordinate system such that 
the $\,\nu$-equations, (\ref{fredh2}),  
become a homogeneous linear system of equations for one unknown function. 
This system, (\ref{nu2}), can be completed into an involutive system. 
\end{theorem}

\begin{remark}\label{two} Note that we want to preserve the 
relationship (\ref{lamdef}), i.e.,
\begin{equation}
\Xi_{i}\,\equiv\,g_{1i}\,-\,\lambda^j_1\,\wg_{ji}\,=\,0.
\nonumber
\end{equation}
Using only the  $\wg$-equations, one computes
\begin{eqnarray}
0\, & = & \,\lambda^j_1\,\{\,\lambda^\ell_1\,
{\partial \wg_{ij}\over\partial x^\ell} 
\,+\,
{\partial  \lambda^\ell_a \over\partial x^i}\cdot\wg_{\ell j}\,+\,
{\partial  \lambda^\ell_a \over\partial x^j}\cdot\wg_{\ell i}\,-\,
{\partial  g_{ij} \over\partial x^1} \,\}\nonumber \\
 & = & \,\lambda^\ell_1\,{\partial \,\Xi_{i}  \over\partial x^\ell}\,+\,
{\partial\lambda^\ell_1\over\partial x^i}\;\Xi_{\ell}\,+\,
{\partial\lambda^\ell_1\over\partial x^i}\,g_{\ell 1}\,+\,
\lambda^\ell_1\,{\partial g_{i1} \over\partial x^\ell} \,-\,
\lambda^j_1\,{\partial g_{ij} \over\partial x^1} \,\nonumber \\
 & = & \,\lambda^\ell_1\,{\partial \,\Xi_{i} \over\partial x^\ell}\,+\,
{\partial\lambda^\ell_1\over\partial x^i}\,\Xi_{\ell}\,
+\,{\partial \hfil \over\partial x^i}\left(\lambda^\ell_1\,g_{\ell 1}\right)\,
-\,\lambda^\ell_1\,\left({\partial g_{\ell 1}\over\partial x^i}\,+\,
{\partial g_{i \ell}\over\partial x^1}\,
-\,{\partial g_{i 1}\over\partial x^\ell}\right) \nonumber 
\end{eqnarray} 
Now invoke the $\lambda$-equation to obtain
\begin{equation}
\lambda^\ell_1\,{\partial \,\Xi_{i}\over\partial x^\ell}\,+\,
{\partial\lambda^\ell_1\over\partial x^i}\;\Xi_{\ell}\,=0.
\nonumber
\end{equation}
We see that equality (\ref{lamdef}) holds locally provided it holds on a 
hypersurface transverse to $\,\lambda^i_1$. See also 
\cite[Proposition 1.4]{AKW}
\end{remark}

Given any non-zero solution of the $\lambda$-equations, 
there is a local coordinate system $\,y^1$, \dots, $y^n$ such that 
\begin{equation}
\lambda^i_1(x)\,{\partial y^j\over\partial x^i}\,=\,\delta^j_1. 
\label{char1}
\end{equation}
Let $G$ and $\widehat G$ represent $g$ and $\wg$ in the $y$-coordinates. 
The equation (\ref{lamdef}) then reads
\begin{equation} 
G_{ij}(y)\,{\partial y^i\over\partial x^1}\,
{\partial y^j\over\partial x^r}\,=\,
\widehat G_{1k}\,{\partial y^k\over\partial x^r}.
\label{lam3}
\end{equation}
The $\wg$-equations read
\begin{equation}
{\partial \widehat G_{ij}\over\partial y^1}\,=\,
{\partial g_{k\ell}\over\partial x^1}\,
{\partial x^k\over\partial y^i}\,
{\partial x^\ell\over\partial y^j}\,,
\label{gh2}
\end{equation}
and the $\wV$-equations become
\begin{equation}
{\partial \widehat V\over\partial y^1}\,=\,
{\partial V\over\partial x^1}\,.
\nonumber
\end{equation}
It is easy to see that the following result holds.
\begin{theorem} 
Given any non-zero solution to the $\lambda$-equation, 
there is a unique solution to the $\wg$- and $\wV$-equations with 
initial data prescribed at $y^1=0$.
\end{theorem}
\begin{remark}\label{three} Note that equation(\ref{lam3}) gives directly
\begin{equation}
\widehat G_{k1}\,=\,G_{ki}\,{\partial y^i\over\partial x^1}\,,
\nonumber
\end{equation}
and so one only needs to solve (\ref{gh2}) for $n\,(n-1)/2$ quantities 
$\,\widehat G_{ij}$, $\;\;2\le i, j\le n$.
\end{remark}

\section{Example 1: Inverted pendulum in a vertical plane} 

As the  first example we consider the inverted pendulum restricted 
to a vertical plane with horizontal and vertical actuation of the base, 
see Figure 1.


\hskip 50bp\psfig{file=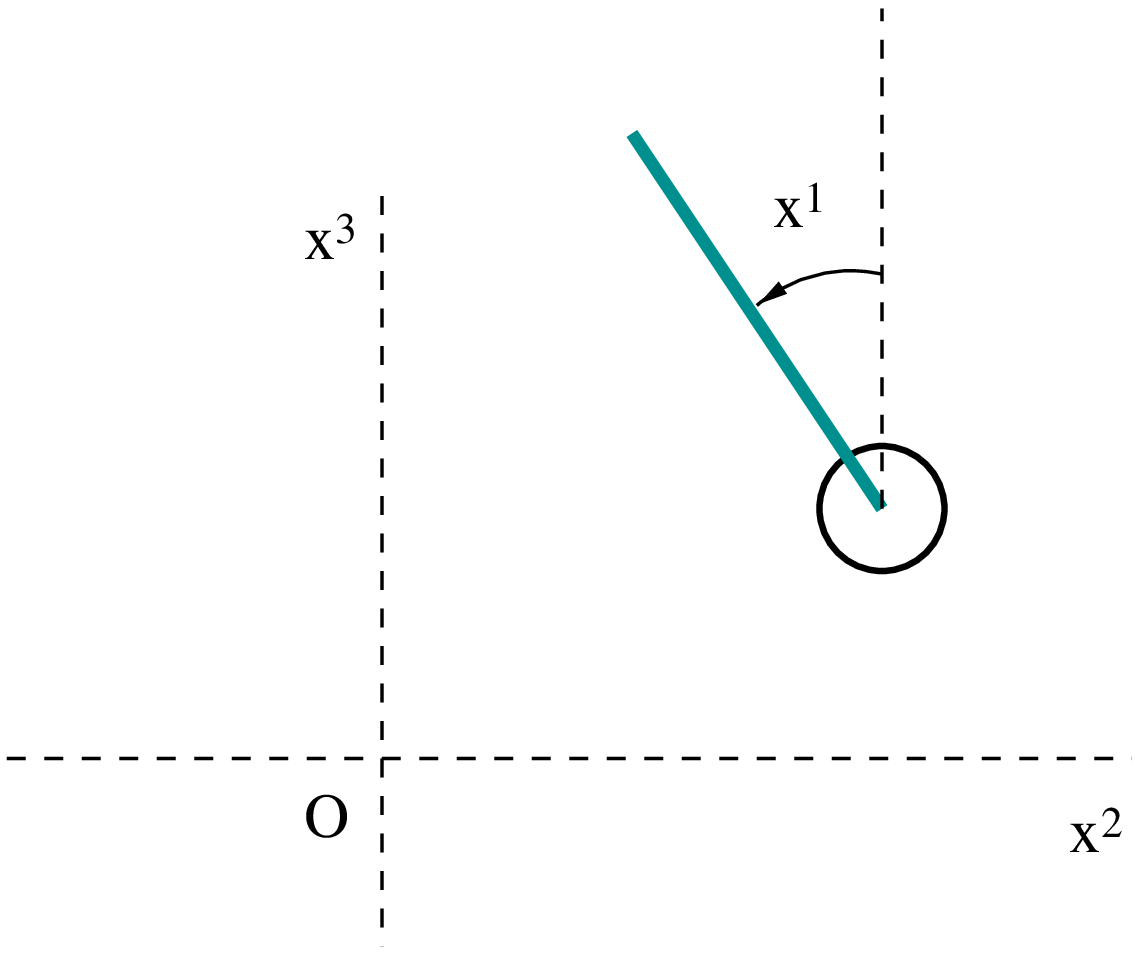,width=3truein,height=2.5truein}%

\centerline{\bf Figure 1\/}

\vfill
\pagebreak

\noindent After rescaling units, the mass matrix and potential energy are given by
\begin{equation}
g\,=\,\begin{pmatrix}
1 & - a\,\cos(x^1) & - a\,\sin(x^1) \\
- a\,\cos(x^1) & 1 & 0 \\
- a\,\sin(x^1) & 0 & 1
\end{pmatrix}
\nonumber
\end{equation}
\begin{equation}
V\,=\,b\,x^3\,+\,\cos(x^1)
\nonumber
\end{equation}
Since only $x^1$ is unactuated, we will simplify notation and use 
$\,\lambda^i\,$ to denote $\,\lambda^i_1$. 
The $\lambda$-equations (\ref{lameq}) are 
$$
\partial_1\,\nu\,=\,2a\,\sin(x^1)\,\lambda^2_1\,-\,2a\,\cos(x^1)\,
\lambda^3_1,\quad
\partial_2\,\nu\,=\,0,\quad \partial_3\,\nu\,=\,0
$$
with 
$\nu\,=\,\lambda^1\,-\,a\,\cos(x^1)\,\lambda^2\,
-\,a\,\sin(x^1)\,\lambda^3$. It is not difficult to see that 
the general solution to these equations is
\begin{eqnarray}
\lambda^1\,& = & \,\nu(x^1)\,+\,\frac12\,\cot(x^1)\,\partial_1\nu(x^1)\,
+\,a {\lambda^3(x^1,x^2,x^3)\over \sin(x^1)} \nonumber\\
\lambda^2\,& = &\,{1\over 2a\,\sin(x^1)}\,\partial_1\nu(x^1)\,
+\,\cot(x^1)\,\lambda^3(x^1,x^2,x^3) \nonumber\\
\,\nonumber \\
\lambda^3\,& = & \,\lambda^3(x^1,x^2,x^3)\,,\nonumber
\end{eqnarray}
where $\,\nu(x^1)$, $\,\lambda^3(x^1,x^2,x^3)\,$ are arbitrary.
In order to obtain a managable explicit solution to the matching 
equations, we will choose 
$$
\nu(x^1)\,=\,a\,\mu_0\,\sin^2(x^1)\,+\,\sigma_0\,-\,a\,\mu_0,
\qquad \lambda^3\,=\,0
$$
with free parameters $\,\sigma_0\,$ and $\,\mu_0$. Then
$$
\lambda^1\,=\,\sigma_0\;,\qquad \lambda^2\,=\,\mu_0\,\cos(x^1).
$$
The coordinates  
$$
y^1\,=\,{1\over\sigma_0}\,x^1\;,
\quad y^2\,=\,x^2\,-\,\mu_0\,\sin(x^1)\;,\quad y^3\,=\,x^3
$$
satisfy (\ref{char1}). Following Remark \ref{three}, we need to solve  
the  $\widehat g$-equations only for $\,\wg_{22}$, $\,\wg_{23}$, 
and $\,\wg_{33}$. These equations are  
$$
{\partial\hfil\over\partial y^1}\,\widehat g_{22}\,=\,
{\partial\hfil\over\partial y^1}\,\widehat g_{23}\,=\,
{\partial\hfil\over\partial y^1}\,\widehat g_{33}\,=\,0.
$$
Clearly,  
\begin{eqnarray}
\widehat g_{22}\,& = & \,\widehat g_{22}(y^2,y^3)\,=\,
\widehat g_{22}(x^2\,-\,\mu_0\,\sin(x^1),\;x^3)\nonumber\\
\widehat g_{23}\,& = & \,\widehat g_{23}(x^2\,-\,\mu_0\,\sin(x^1), \;x^3)
\nonumber\\
\widehat g_{33}\, & = & \,\widehat g_{33}(x^2\,-\,\mu_0\,\sin(x^1), \;x^3)
\nonumber
\end{eqnarray}
From $\,g\,=\,\widehat g\,\lambda$, we obtain 
the rest of $\wg_{ij}$ as:  
\begin{eqnarray}
\widehat g_{11}\, & = & \,{1\over\sigma_0}\,
+\,{a\mu_0\over\sigma_0^2}\,\cos^2(x^1)\,
+\,{\mu_0^2\over\sigma_0^2}\,\cos^2(x^1)\;\widehat g_{22}\nonumber\\ 
\widehat g_{12}\, & = & \,-\,{a\over\sigma_0}\,
-\,{\mu_0\over\sigma_0}\,\cos(x^1)\;\widehat g_{22}\nonumber\\ 
\widehat g_{13}\, & = & \,-\,{a\over\sigma_0}\,\sin(x^1)\,
-\,{\mu_0\over\sigma_0}\,\cos(x^1)\;\widehat g_{23}\nonumber
\end{eqnarray}
The $\widehat V$-equation yields
$$
\widehat V\,=\,{1\over\sigma_0}\,\cos(x^1)\,+\,w(y^2, y^3)\,.
$$
The $\widehat C$-equation reads\ \ \ \ 
$\lambda^j\,\widehat C_{j }\,=\,0$. One solution is
$$
\widehat C\,=\,-\,\sigma_0\,R(x)\;
\begin{pmatrix}  {\mu_0^2\over\sigma_0^2}\,\cos^2(x^1)\, &
-\,{\mu_0\over\sigma_0}\,\cos(x^1) & -\,{\mu_0\over\sigma_0}\,\cos(x^1) \\
-\,{\mu_0\over\sigma_0}\,\cos(x^1) & 1 & 1\\
 -\,{\mu_0\over\sigma_0}\,\cos(x^1) & 1 & 1
\end{pmatrix}\cdot 
\begin{pmatrix}  \dot x^1 \\ \dot x^2 \\ \dot x^3
\end{pmatrix} 
$$
The resulting control law can be obtained explicitly from equation (\ref{control}). 
The expression is too long to be included in this paper.

\begin{proposition} 
If the functions 
$\widehat g_{22}(y^2, y^3)$, $\widehat g_{23}(y^2, y^3)$,
$\widehat g_{33}(y^2, y^3)$, $w(y^2,y^3)$,  and $R(x)$, 
and the 
parameters $\mu_0$ and $\sigma_0$   
are chosen so that
\begin{equation}
\begin{split}
\widehat g_{22}(0)\,>\,0\;,\quad \widehat g_{23}(0)\,=\,0\,,
\quad \widehat g_{33}(0)\,=\,1, \\
\partial_{y^2}^2 w(0)\,>\,0, \quad \partial_{y^2}\partial_{y^3} w(0)\,=\,0,
\quad \partial_{y^3}^2 w(0)\,>\,0,\qquad R(0)\,>\,0,\\
\sigma_0\,<\,0,\quad \widehat g_{22}(0)\,\mu_0^2\,+\,a\,\mu_0\,+\,
\sigma_0\,>\,0\;,
\quad  \widehat g_{22}(0)\,(a\,\mu_0\,-\,\sigma_0)\,+\,a^2\,<\,0\;,
\end{split}
\nonumber
\end{equation}
then $\,x\,=\,\dot x\,=\,0\,$ is a 
locally 
asymptoitcally stable equilibrium of the closed loop system.
\end{proposition}

\section{Example 2: Inverted pendulum cart on a seesaw }

In the previous example the kernel of the matrix 
$\,A^\star\,$, (\ref{nualg}), 
was two-dimensional. Generically, for systems with one 
unactuated degree of freedom the dimension of the kernel 
will be $1$. The following example illustrates this situation. 
The inverted pendulum cart on a seesaw is shown in Figure 2. 
There are several interesting ways to actuate this system. 
We will consider the case with actuated cart and pendulum, 
and unactuated seesaw.


\psfig{file=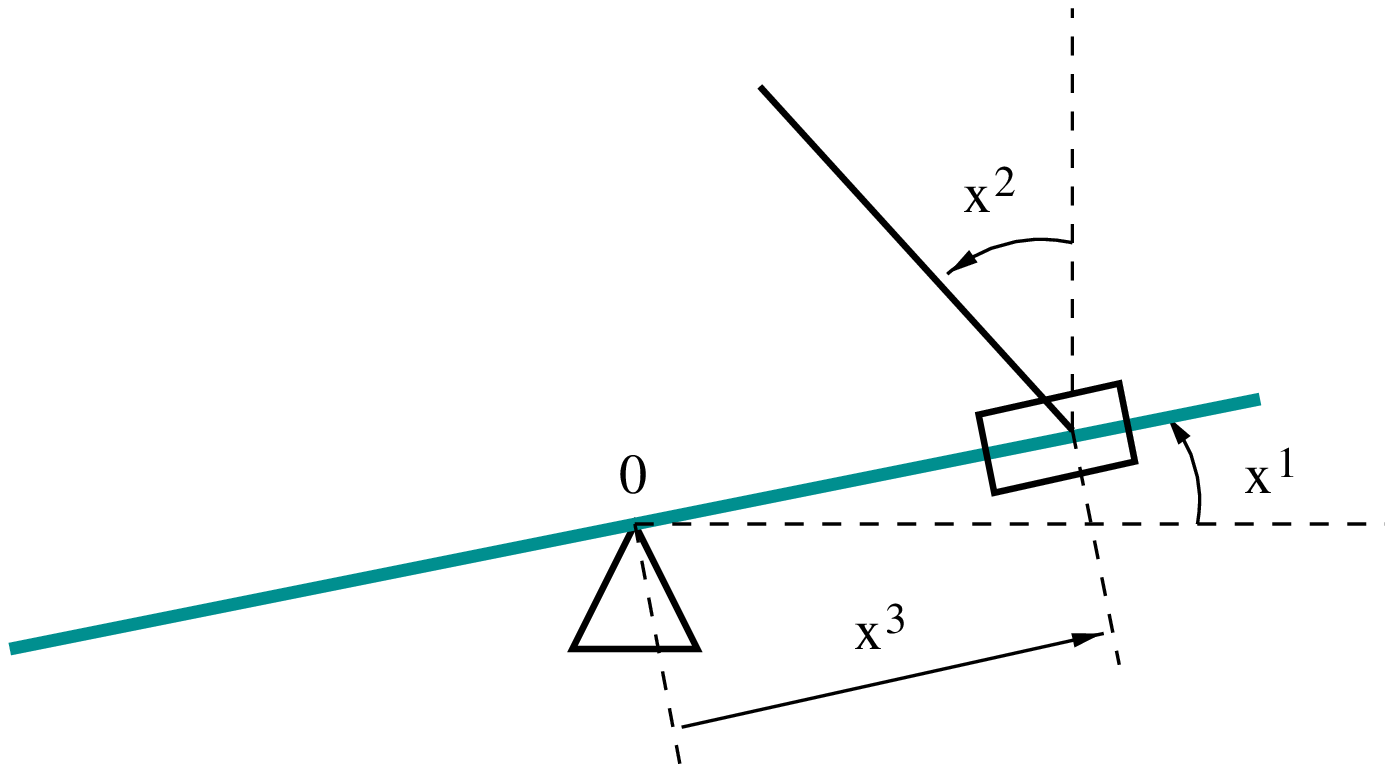,width=4.5truein,height=2truein}%

\centerline{\bf Figure 2\/}

\noindent The rescaled mass matrix and potential energy of the system are given by
\begin{equation}
g\,=\,
\begin{pmatrix} b\,+\,(x^3)^2 & a\,x^3\,\sin(x^1-x^2) & 0 \\
a\,x^3\,\sin(x^1-x^2) & 1 & -\,a\,\cos(x^1-x^2) \\
0 & -\,a\,\cos(x^1-x^2) & 1
\end{pmatrix}
\nonumber
\end{equation}
and 
\begin{equation}
V\,=\,x^3\,\sin(x^2)\,+\,a\,\cos(x^1)\,. 
\nonumber
\end{equation}
The theory in section 3 was presented with special coordinates 
so that $\,g_{11}=1$.  However, in practice this is not necessary.  

As before, we write $\,\lambda^i$ for $\,\lambda^i_1$ 
and  $\,\nu\,$ for $\,g_{1j}\,\lambda^j$.  
The $\lambda$-equations are
\begin{eqnarray}
\partial_1\,\nu\, & = & \,2a\,x^3\,\cos(x^1-x^2)\,\lambda^2\,
-\,2\,x^3\,\lambda^3 \nonumber\\ 
\partial_2\,\nu\, & = & \,0 \nonumber\\
\partial_3\,\nu\, & = & \,2a\,\sin(x^1-x^2)\,\lambda^2\,
+\,2\,x^3\,\lambda^1 \nonumber
\end{eqnarray}
Hence, $\,\nu\,=\,\nu(x^1,x^3)$. 
Plug in  
$\,\lambda^1 \,=\,(\nu-g_{12}\,\lambda^2-g_{13}\,\lambda^3)/g_{11}$ and 
solve for $\,\lambda^2\,$ and $\,\lambda^3$:
\begin{eqnarray}
\lambda^1\, & = & \,{1\over 2b}\,\big(\,2\nu\,-\,x^3\partial_3\,\nu\,)
\nonumber\\ 
\lambda^2\, & = & \,{1\over 2ab\,\sin(x^1-x^2)}\,
\bigg(-\,2\,x^3\,\nu\,+(b\,+\,(x^3)^2)\,\partial_3\,\nu\,\bigg)\nonumber\\ 
\lambda^3\, & = & \,{1\over 2b\,x^3\,\sin(x^1-x^2)}\,
\bigg(-\,2(x^3)^2\,\cos(x^1-x^2)\,\nu\, \nonumber\\
& + &\,x^3\,(b+(x^3)^2)\,\cos(x^1-x^2)\,\partial_3\nu\,
-\,b\,\sin(x^1-x^2)\,\partial_2\nu \bigg)\nonumber
\end{eqnarray}

Notice, that $\,\lambda^2\,$ and $\,\lambda^3\,$ blow up as $\,x\,$ 
approaches $0$ unless $\,\nu=\kappa\,(b+(x^3)^2)$. Since 
$\,g\,=\,\wg\,\lambda$, one must have $\,\det\wg\to 0$ as $\,x\to 0$, 
i.e., $\,\wg$ degenerates at $x=0$.  
This means that 
$\,\widehat H(x, \dot x)\,=\,\frac12 \,\wg_{ij}\,\dot x^i\,\dot x^j\,+\,\widehat V\,$ 
cannot serve as a Lyapunov function unless $\,\nu=\kappa\,(b+(x^3)^2)$. 
This $\nu$ corresponds exactly to the basic solutions of the matching 
equations from Remark \ref{remarkone}. This illustrates the following 
general principle.

\begin{remark}\label{princ} If $\,(x_0,\,0)\,$ is the desired 
equilibrium of a system and 
\begin{equation}
\hbox{rank}\,A^\star(x_0)\,<\,\limsup_{x\to x_0} 
\hbox{rank}\,A^\star(x)\,,
\label{rankdrop}
\end{equation}
then only  basic solutions of the matching equations should be tested 
to produce a stabilizing control law from (\ref{control}).
\end{remark}

\section{Example 3: inverted pendulum cart on a roller coaster}

Consider a cart with inverted pendulum on a roller coaster.
Special cases of this mechanical system include the inverted
pendulum on a rotor arm, the inverted pendulum on a verticle disk, and 
the inverted pendulum cart on an incline. 
By assuming that the size of the base of the cart is relatively small 
we may neglect the inertia of the base of the cart. 
It is therefore sufficient to model the cart with one point mass 
 for the base, 
and one point mass a fixed distance 
away for the pendulum. The pendulum joint will be unactuated.

The configuration of the system may be described by a position and 
an angle. Assume that the shape of the roller coaster is given 
as a curve $\,x(s)\,$ in $\,\mathbb R^3\,$ parametrized by 
arc length, $\,s$, from a fixed point. 
Assume that the pendulum is always in the plane 
containing the tangent vector, $\,\tau(s)$, and the vertical direction, 
$\,e_3$. Let $\,\phi\,$ be the angle between the pendulum 
and $\,e_3$. By rescaling mass, length, and time, we will write 
\begin{eqnarray}
g\, & = & \begin{pmatrix} 1 & 
b\,\sin(\alpha-\phi)\nonumber\\
b\,\sin(\alpha-\phi) & 
\left(1+k(s)^2\,\sin^2\phi\,
\big({\sin^2\alpha-n_3^2\over \sin^4\alpha}\big)\right)
\end{pmatrix} \nonumber\\
V\, & = & \,a\,x^3\,+\,\cos\phi,\nonumber
\end{eqnarray}
where $a$ and $b$ are positive parameters, $\,0<b<1$, 
and $x^3$ is the vertical 
component of $x(s)$ . 
The (unit) tangent vector to the curve is 
$\,\tau(s)={x^\prime(s)\over |x^\prime(s)|}$, where $\,{}^\prime\,$ 
stands for the derivative with respect to $s$. 
The curvature of the curve 
is $\,k(s)=|\tau^\prime(s)|$. Denote by $\,n(s)\,$ the principal 
normal to the curve. Recall that $\,\tau^\prime(s)=k(s)\,n(s)$. 
In the formula above $\,n_3\,$ is the vertical component of the 
principal normal, and $\alpha$ is the angle between $\,\tau\,$ 
and the vertical direction. 
Index $1$ corresponds to $\phi$, index $2$ corresponds to $s$.
The unactuated degree of freedom corresponds to the $\phi$ variable. 
The $\lambda$-equations (\ref{nu}) then read as follows:
\begin{equation}
\begin{aligned}
\partial_1\,\nu\, & =  -\,2\,b\,\cos(\alpha - \phi)\;\lambda^2_1 \\ 
\partial_2\,\nu\, & =  \;k(s)^2\,\sin(2\,\phi)\,
\big({\sin^2\alpha-n_3^2\over \sin^4\alpha}\big)\;\lambda^2_1\,.
\end{aligned}
\label{dnu}
\end{equation}
The orthogonality equation (\ref{nu2}) then, obviously, is
\begin{equation}
k(s)^2\,\sin(2\,\phi)\,
\big({\sin^2\alpha-n_3^2\over \sin^4\alpha}\big)\;
{\partial\nu\over \partial\phi}\,+\,
2\,b\,\cos(\alpha - \phi)\;{\partial\nu\over \partial s}\,=\,0\,.
\label{rol1}
\end{equation}
It is not clear if all solutions to this equation can be written 
explicitly for 
a  general curve. We consider here two particular cases when this 
is possible. 
The first case is when $\sin^2\alpha=n_3^2$. 
This occurs exactly when the roller coaster lies in one vertical plane. 
The second case is when $\,\alpha(s)\,$ is constant. This occurs when the 
track is constantly inclined. 

\subsection{Case 1: $\sin^2\alpha\,=\,n_3^2$}
Note that this case 
includes the interesting examples of an inverted pendulum on a 
vertical disk and an inverted pendulum cart on an incline. 

As is readily seen from (\ref{rol1}), the general solution of 
(\ref{rol1}) in this case is $\,\nu=\nu(\phi)$, an arbitrary  function. 
Then 
\begin{equation}
\lambda^2_1\,=\,-{1\over 2\,b\,\cos(\alpha - \phi)}\; 
{\partial\nu\over \partial\phi}
\nonumber
\end{equation}
and
\begin{equation}
\lambda^1_1\,=\,\nu(\phi)\,
+\,\frac12\,\tan(\alpha - \phi)\; {\partial\nu\over \partial\phi}.
\nonumber
\end{equation}
This is a general solution of the $\lambda$-equation. From here 
one must solve the $\wg$- and $\wV$-equations. For special choices 
of $\,\alpha(s)\,$ and/or $\,\nu(\phi)\,$ these equations have 
explicit closed form solutions. 

\subsection{Case 2: $\alpha(s)\,=\,\alpha_0$}

Examples with $\,\alpha(s)\,$ constant include an inverted pendulum cart 
traveling on any path in a horizontal plane, 
a cart on a vertically oriented circular helix, or 
any constantly 
inclined track.

Since ${d\alpha\over ds}=-k(s){n_3(s)\over \sin(\alpha)}$, 
if $\,\alpha(s) =\alpha_0$, then we have $k(s)n_3(s)=0$. 
To solve equation (\ref{rol1}), we introduce new coordinates 
$$
\aligned 
z^1\,& =\,\beta(s)\\ 
z^2\,& =\,
\beta(s)+\cos(\alpha_0)\ln |\csc\phi+\cot\phi|-
\sin(\alpha_0)\ln |\sec\phi+\tan\phi |,
\endaligned
$$
where 
$$
\beta(s)=\int_0^s {k^2(p)\over b\sin^2(\alpha_0)}\,dp\,.
$$
Equation (\ref{rol1}) is equivalent to the fact that $\,\nu\,$ 
is an arbitrary 
function of $\,z^2$. Thus, from (\ref{dnu}) and the definition of $\,\nu\,$ 
we find
$$
\lambda^1_1\,=\,\nu(z^2)
-{\sin(\alpha_0-\phi)\over \sin(2\phi)}\,{d\,\nu\over d\,z^2}\,;\;\qquad 
\lambda^2_1\,=\,{1\over b\,\sin(2\phi)}\,{d\,\nu\over d\,z^2}\,.
$$

\subsection{End of the roller coaster example}

From the computations in case 1 and case 2 one can see that 
the general solution to the matching equations for the cart 
on a roller coaster will be fairly complicated. However, we can show that 
any {\it linear\/} control law is the first order germ (linearization) 
of some matching control law. In fact, the only requirement for this 
is that there is no rank drop at the equilibrium, i.e.,
\begin{equation}
\hbox{rank}\,A^\star(x_0)\,=\,\limsup_{x\to x_0} 
\hbox{rank}\,A^\star(x)\,.
\label{norankdrop}
\end{equation}
We assume that the dissipative term at the equilibrium satisfies 
the following natural assumptions:
\begin{equation}
C_\ell(x_0, 0)\,=\,0,\qquad 
{\partial\hfil\over\partial x^i}\,\bigg|_{(x_0, 0)}\;C_\ell\,=\,0\,.
\nonumber
\end{equation}
\begin{lemma} If condition (\ref{norankdrop}) is satisfied for a two degree 
of freedom system, then the first order germs of matching control laws 
 at $\,(x_0, 0)\,$ exhaust all linear control laws for which 
the closed loop system has an equilibrium at $\,(x_0, 0)\,$.
\end{lemma}
{\em Proof}. Given a linear control input
$$
u^{\hbox{lin}}_j\,=\,v_j\,+\,a_{ij}\,(x^i-x^i_0)\,+\,b_{ij}\,\dot x^i
$$
with $\,u^{\hbox{lin}}_1\,=\,0$, we will find a matching 
control law with the same 
germ. From the general expression (\ref{control}) 
for the matching control law, we see that the first order germ is 
$$
u^{\hbox{germ}}_j\,=\,(V_j\,- \,g_{j \ell}\,
\wg^{\ell r}\,\widehat V_r)\,   
+\,\big(\, V_{j r}\,-\,g_{j \ell}\,
\wg^{\ell i}\,\widehat V _{ir}\,\big)\,(x^r-x_0^r)\,
+\, 
\, \big(\, C_{j i}\,-\,\,g_{j \ell}\,
\wg^{\ell r}\,\widehat C_{r i}\,\big)\,\dot x^i,
$$
where 
$$ 
V_{j}\, 
= \,{\partial V\over\partial x^j}\bigg|_{x_0}\,,\quad  
V_{j r}\, 
= \,{\partial^2 V\over\partial x^j\partial x^r}\bigg|_{x_0}\;,\quad
C_{j i}\,=\,{\partial C_j\over\partial \dot x^i}\bigg|_{(x_0, 0)}\,, 
$$
and $\,\wV_j$, $\,\widehat V_{jr}$, and 
$\,\widehat C_{ji}\,$  are defined similarly. Equating like terms 
gives
\begin{equation}
\wV_{\ell j}\,=\,\wg_{\ell i}\,g^{ir}\,(V_{rj} - a_{rj})\,,\quad
\widehat C_{\ell j}\,=\,\wg_{\ell i}\,g^{ir}\,(C_{rj} - b_{rj})\,,\quad
\wV_\ell\,=\,\wg_{\ell i}\,g^{ir}\,(V_{r} - v_r)=0\,.
\label{compar}
\end{equation}
One can see that $\,\wV_j$, $\,\widehat V_{jr}$, and 
$\,\widehat C_{ji}\,$ are specified once $\,\wg_{ij}(x_0)\,$ is known. 
Moreover, there exists a non-degenerate symmetric $\,\wg_{ij}(x_0)\,$ 
such that the resulting $\,\wV_{\ell j}\,$ will be symmetric, see 
\cite[Lemma 1]{BBeam2}.  
To conclude the argument, we now show that any 
non-degenerate symmetric $\,\wg_{ij}(x_0)\,$ arises as a zero order 
germ of a solution to $\,\wg$-equation. Also, any $\,\wV_\ell\,$, 
$\,\wV_{\ell j}\,$  
satisfying $\,\wV_\ell\,=\,\wg_{\ell i}\,g^{ir}\,(V_{r} - v_r)\,$ 
and $\,\wV_{\ell j}\,=\,\wg_{\ell i}\,g^{ir}\,(V_{rj} - a_{rj})\,$ 
arises as a solution to the $\,\wV$-equation. 

Given a non-degenerate $\,\wg_{ij}(x_0)$, define the non-degenerate 
$\,\lambda^j_i(x_0)=g_{ik}(x_0)\,\wg^{kj}(x_0)$. Set 
$\,\nu_0\,=\,g_{11}(x_0)\,\lambda^1_1(x_0)+g_{12}(x_0)\,\lambda^2_1(x_0)$. 
The $\lambda$-equations 
in this case are
\begin{equation}
\begin{aligned}
\partial_1 \nu 
- 2\,[1\,1,\,2]\,{1\over g_{11}}\,\nu\; & = & \;2\,[1\,1,\,2]\,
\lambda^2_1 \\
\partial_2 \nu 
- 2\,[1\,2,\,1]\,{1\over g_{11}}\,\nu\; & = &\; 2\,[2\,1,\,2]\,
\lambda^2_1 
\end{aligned}
\label{nuequa}
\end{equation}
By the rank condition, we know that the 
rank of $\,A^\star\,$ in the neighborhood of $x_0$ is either identically 
$0$ or identically $1$. If this rank is $0$, then $\,\nu\,$ can be any 
constant times $\,g_{11}$. We simply choose 
$\,\nu(x)\,=\,(\nu_0/g_{11}(x_0))\,g_{11}(x)$.  
Any solution to the algebraic equation 
$\,\nu(x)\,=\,g_{11}(x)\,\lambda^1_1(x)+g_{12}(x)\,\lambda^2_1(x)$ is 
a solution to the $\lambda$-equation. If the rank of $\,A^\star\,$ is $1$, 
then the orthogonality condition, (\ref{fredh}), is
\begin{equation}
 [2\,1,\,2]\,\partial_1 \nu \,-\,[1\,1,\,2]\,\partial_2 \nu 
+ 2\,\left( [1\,1,\,2]\,[1\,2,\,1]\,
- \,[2\,1,\,2]\,[1\,1,\,2]\,\right)\,{1\over g_{11}}\,\nu \,=\,0\,.
\label{ortho2}
\end{equation}  
At $\,x=x_0$, either 
$\,\partial_1\,$ or $\,\partial_2\,$ is not parallel to the vector 
$\,[2\,1,\,2]\,\partial_1\,-\,[1\,1,\,2]\,\partial_2$ . Assume it is 
$\,\partial_1\,$. Then initial conditions to equation (\ref{ortho2}) 
can be specified along the line $\,x^2=x^2_0$. In particular, 
we can choose the initial values so that 
\begin{equation}\label{bound}
\left( \partial_1 \nu - 2\,[1\,1,\,2]\,{1\over g_{11}}\,\nu \right)\,
\bigg|_{x=x_0}\;=\,2\,[1\,1,\,2]\,\bigg|_{x=x_0}\;\lambda^2_1(x_0)\,,\qquad
\nu(x_0)\,=\,\nu_0\,.
\end{equation}
The second equation in (\ref{nuequa}),
$$
\left( \partial_2 \nu - 2\,[1\,2,\,1]\,{1\over g_{11}}\,\nu \right)\,
\bigg|_{x=x_0}\;=\,2\,[2\,1,\,2]\,\bigg|_{x=x_0}\;\lambda^2_1(x_0)\,,
$$
will be satisfied automatically since 
the rank of $\,A^\star\,$ is $1$. Let $\,\nu(x)\,$ be a solution 
of (\ref{ortho2}) with initial condition, $\,\nu(x^1,0)\,$,  
satisfying (\ref{bound}). Now, one solves (\ref{nuequa}) for 
$\,\lambda^2_1(x)$ 
and then 
$\,\nu(x)\,=\,g_{11}(x)\,\lambda^1_1(x)+g_{12}(x)\,\lambda^2_1(x)\,$ 
for $\,\lambda^1_1(x)$.

\smallskip
Now that the $\lambda$-equations are solved, we turn to $\wg$-equations. 
These equations take the form
$$
{\partial\hfil\over\partial y^1}\,\wg\,+\,R\,\wg\,=\,S\,,
$$
where $\,{\partial\hfil\over\partial y^1}\,=\,
\lambda^1_1\,\partial_1\,+\,\lambda^2_1\,\partial_2$. The initial 
conditions can be set on any line transverse to 
$\,{\partial\hfil\over\partial y^1}$, in particular, along the line 
$\,\lambda^1_1(x_0)\,(x^1-x^1_0)\,+\,\lambda^2_1(x_0)\,(x^2-x^2_0)=0$. 
On this line set $\,\lambda^i_2(x)\,=\,\lambda^i_2(x_0)$, and 
$\,\wg(x)\,=\,g(x)\cdot(\lambda(x_0))^{-1}$. The solution to 
equation (\ref{gheq}) with this initial data then has the desired value at  
$\,x=x_0$.

\smallskip
It remains to show that the $\,\wV$-equation has a solution such that 
$\,\wV_\ell\,=\,0\,$ and 
\begin{equation}
\lambda^\ell_r(x_0)\,\wV_{\ell j}\,=\,V_{rj} - a_{rj}\,,
\label{vlj}
\end{equation}
where $\,\lambda^\ell_r(x_0)=g_{rk}(x_0)\,\wg^{k\ell}(x_0)$. Since 
$\,\lambda^\ell_r(x_0)\,$ is non-degenerate, either 
$\,\lambda^1_1(x_0)\ne 0\,$ or $\,\lambda^2_1(x_0)\ne 0$. Consider the case 
with $\,\lambda^2_1(x_0)\ne 0$. In this case the line $\,x^2=x^2_0\,$ is 
non-characteristic for the $\,\wV$-equation 
\begin{equation}
\lambda^1_1\,\partial_1\,\wV\,+\,\lambda^2_1\,\partial_2\,\wV\,
=\,\partial_1 V\,.
\label{vhat2} 
\end{equation}
Pick the initial value, $\,\wV\big|_{x^2=x^2_0}\,$, so that 
$$
\,\wV_1\,=\,0, \quad 
\wV_{11}\,
=\,{\lambda^2_2\,V_{11}-\lambda^2_1\,(V_{12}-a_{12})\over
\lambda^1_1\,\lambda^2_2\,-\,\lambda^1_2\,\lambda^2_1}\,\bigg|_{x=x_0}
$$ 
and solve equation (\ref{vhat2}). Since $\,x=x_0\,$ is an equilibrium, 
$\,V_1=0$ and $\,(V_2 - v_2)=0$. From the differential equation,
(\ref{vhat2}), we see that $\,\wV_2=0$. Differentiating equation 
(\ref{vhat2}) with respect to $\,x^1\,$ and $\,x^2$,  we see that 
$\,W_{ij}\,=\,\wV_{ij}\,$ satisfies
\begin{equation}
\begin{aligned}
\lambda^1_1(x_0)\,W_{11}\,
+\,\lambda^2_1(x_0)\,W_{12}\,
& = & \,V_{11}\,\\ 
\lambda^1_1(x_0)\,W_{12}\,
+\,\lambda^2_1(x_0)\,W_{22}\,
& = & \; V_{12}\,.
\end{aligned}
\label{vhat3} 
\end{equation}
By construction, 
\begin{equation}
\lambda^1_2(x_0)\,W_{11}\,
+\,\lambda^2_2(x_0)\,W_{12}\;
 =  \, V_{12} - a_{12}\;.
\label{vhat4}
\end{equation}
Notice that (\ref{compar}) implies that 
$\,W_{\ell j}\,=\,\wg_{\ell i}\,g^{ir}\,(V_{rj} - a_{rj})$ also satisfy 
equations (\ref{vhat3}) and (\ref{vhat4}). 
Since the solution to the algebraic system (\ref{vhat3}), (\ref{vhat4}) 
is unique, we conclude that equation (\ref{vlj}) is valid, as required.

\section{Example 4: A double pendulum on a wheel}

Our next example is the system with two unactuated degrees of freedom 
depicted in Figure 3. Only joint {\bf A\/} is actuated.

After rescaling, the entries $\,g_{ij}\,$ of the mass matrix are
\begin{equation}
g_{ij}\,=\,m_{ij}\,\cos(x^i-x^j)
\nonumber 
\end{equation} 
and the potential energy is 
\begin{equation}
V\,=\, a_1\,\cos(x^1)\,+\,a_2\,\cos(x^2)\,+\,a_3\,\cos(x^3)\,. 
\nonumber 
\end{equation} 
The parameters $\,m_{ij}=m_{ji}\,$ and $\,a_{j}$ are positive.

\bigskip

\hskip 100bp\psfig{file=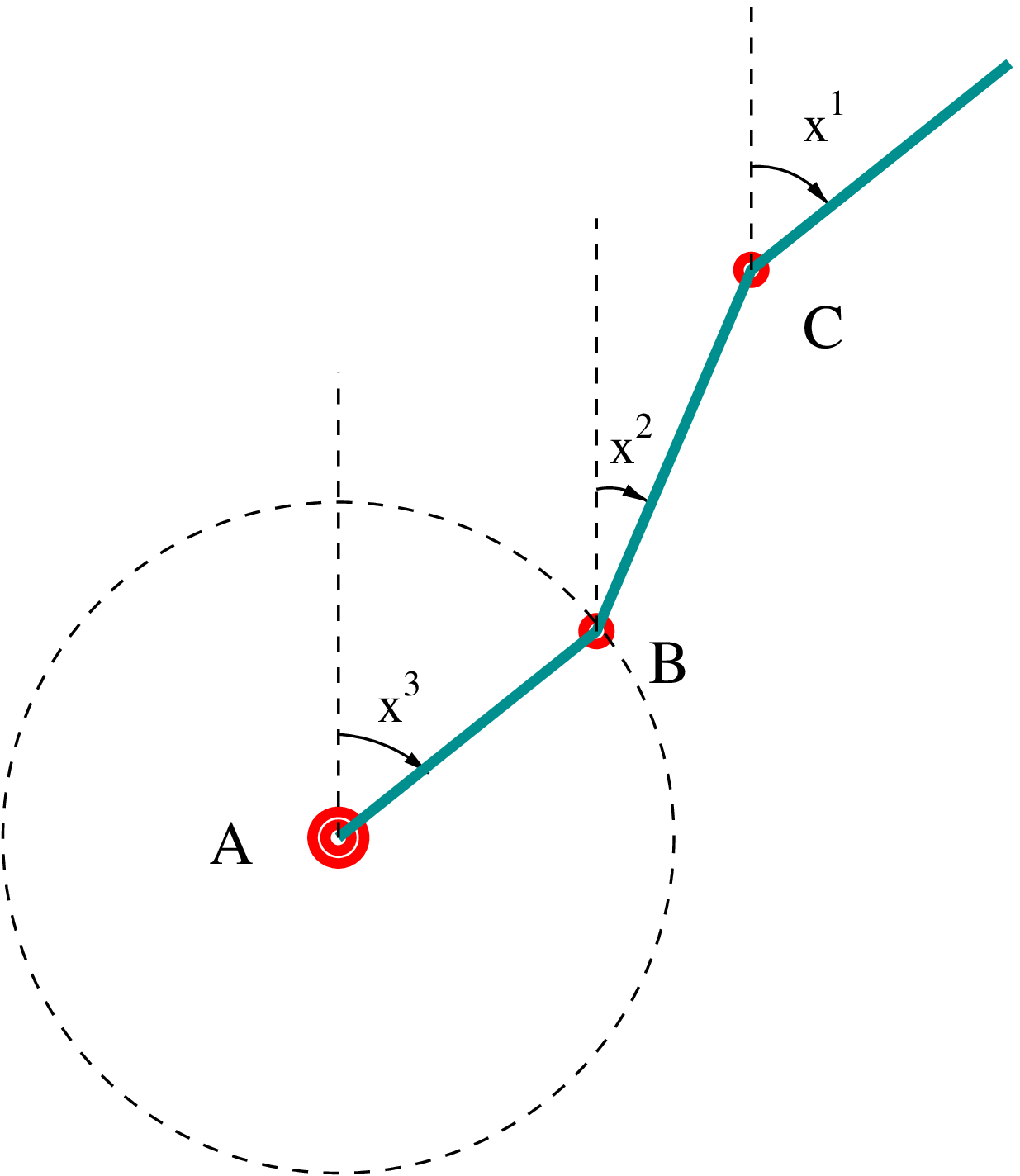,width=2.4truein,height=2.6truein}%

\centerline{\bf Figure 3\/}

\noindent There are six unknown $\,\lambda^i_a$. 
Define $\,\nu_{ab}=g_{ai}\,\lambda^i_b$. Note that we must have 
$\,\nu_{12}=\nu_{21}$. The computations in this section were 
performed using 
Maple. 
The $\lambda$-equations (\ref{lameq}) are
$$
\aligned 
\partial_1 \,\nu_{11}\, & =  \,-\,2\,m_{12}\,\sin(x_1-x_2)\,\lambda^2_1\,
-\,2\,m_{13}\,\sin(x_1-x_2)\,\lambda^3_1  \\ 
 \partial_2 \,\nu_{11}\, & = \,0\,\\ 
 \partial_3 \,\nu_{11}\, & =  \,0 \\
\; \\ 
\partial_1\,\nu_{22}\,  & =   \,0 \\
 \partial_2 \,\nu_{22}\, & =  +\,2\,m_{12}\,\sin(x_1-x_2)\,\lambda^1_2\,
-\,2\,m_{23}\,\sin(x_2-x_3)\,\lambda^3_2   \\ 
\partial_3 \,\nu_{22}\, & =  \,0 \\
\; \\
\partial_1 \,\nu_{12}\, & =   -\,m_{12}\,\sin(x_1-x_2)\,\lambda^2_2\,
-\,m_{13}\,\sin(x_1-x_3)\,\lambda^3_2  \\ 
\partial_2 \,\nu_{12}\, & =  +\,m_{12}\,\sin(x_1-x_2)\,\lambda^1_1\,
-\,m_{23}\,\sin(x_2-x_3)\,\lambda^3_1  \\ 
\partial_3 \,\nu_{12}\, & =  \,0  
\endaligned
$$    
The second step is to express  
$\lambda^1_1$, $\lambda^1_2$, $\lambda^2_1$, and 
$\lambda^2_2$, in terms of $\,\nu_{11}$, 
$\,\nu_{12}$, and $\,\nu_{22}$, 

After substitution  into the above equations  we obtain
\begin{equation}\label{compat}
\aligned
\partial_1 \,\nu_{11}\,&=\,D_{1,\,1}^1\,\nu_{11}\,+\,
D_{1,\,1}^2\,\nu_{12}\,+\,
B_{1,\,1}^1\,\lambda^3_1  \\ 
\partial_2 \,\nu_{11}\,&=\,0\,\\ 
\partial_3 \,\nu_{11}\,&=\,0 \\ 
\partial_1 \,\nu_{22}\,&=\,0 \\
\partial_2 \,\nu_{22}\,&=\,D_{3,\,2}^2\,\nu_{12}\,+\,
D_{3,\,2}^3\,\nu_{22}\,+\,
B_{3,\,2}^2\,\lambda^3_2  \\ 
\partial_3 \,\nu_{22}\,&=\,0 \\
\partial_1 \,\nu_{12}\,&=\,D_{2,\,1}^2\,\nu_{12}\,+\,
D_{2,\,1}^3\,\nu_{22}\,+\,
B_{2,1}^2\,\lambda^3_2 \\ 
\partial_2 \,\nu_{12}\,&=\,D_{2,2}^1\,\nu_{11}\,+\,
D_{2,2}^2\;\nu_{12}\,+\,
B_{2,2}^1\,\lambda^3_1  \\ 
\partial_3 \,\nu_{12}\,&=\,0  
\endaligned    
\end{equation}
Here the $\,D_{i,\,j}^k\,$ and $\,B_{i,\,j}^k\,$ are explicit expressions 
involving $x$. 
In Section 2 we described a general procedure to obtain compatibility 
conditions 
for this system. In this particular case, however, we use a different 
tactic: we 
compute and compare the mixed derivatives of $\,\nu_{ab}$. The first 
set of equations 
we obtain is
$$
\aligned
\partial_3\partial_1\,\nu_{12}\,&=\,K_{11}\,\partial_3\lambda^3_2\,
+\,K_{12}\,\lambda^3_2\,=\,0 \\
\partial_3\partial_2\,\nu_{22}\,&=\,K_{21}\,\partial_3\lambda^3_2\,
+\,K_{22}\,\lambda^3_2\,=\,0 
\endaligned
$$
Direct computation shows that  $\,\det(K_{ij})\,\ne \,0$.  Hence,  
$\,\lambda^3_2\,=\,0$. 
Similarly,
$$
\aligned
\partial_3\partial_1\,\nu_{11}\,&=\,L_{11}\,\partial_3\lambda^3_1\,
+\,L_{12}\,\lambda^3_1\,=\,0 \\
\partial_3\partial_2\,\nu_{12}\,&=\,L_{21}\,\partial_3\lambda^3_1\,
+\,L_{22}\,\lambda^3_1\,=\,0 
\endaligned
$$
and $\,\det(L_{ij})\,\ne \,0$. Hence,  
$\, \lambda^3_1\,=\,0$. 

Next, we substitute 
$\,\lambda^3_1\,=\,\lambda^3_2\,=\,0\;$ into (\ref{compat})  
and solve for $\,\nu_{11}$, $\,\nu_{12}$, $\,\nu_{22}$. This gives
$$
\nu_{11}\,=\,\hbox{const},\qquad
\nu_{22}\,=\,{m_{22}\over m_{11}}\;\nu_{11}\;,\qquad 
\nu_{12}\,=\,{m_{12}\over m_{11}}\,\cos(x^2-x^1)\;\nu_{11}
$$
Returning to $\lambda$-equations 
we see that 
$$
\lambda^1_1\,=\,\lambda^2_2\,=\,{1\over m_{11}}\;\nu_{11}\,,\qquad
\lambda^1_2\,=\,\lambda^2_1\,=\,0\,.
$$
Our computation shows that the only solutions of the matching equations are 
the basic solutions defined in Remark \ref{remarkone}.

\hyphenation{Sprin-g-er-Ver-lag}

\end{document}